\documentclass[12pt]{article}
\usepackage{
amsfonts,amsmath,amssymb, 
geometry,esint}
\DeclareMathOperator\artanh{artanh}
\catcode`\@=11

\@addtoreset{equation}{section}

\usepackage{color}

\def\R{\mathbb{R}}
\def\C{\mathbb{C}}
\def\H{\mathbb{H}}
\def\ms{~\!}
\def\deg{\rm deg}
\def\S{\mathbb{S}}
\def\f{\varphi}
\def\eps{\varepsilon}
\def\beps{\overline\varepsilon}
\def\div{{\rm div}}

\def\irn{\int\limits_{\S^1}\!\!}
\def\fintm{\fint\limits_{\S^1}}

\def\proof{\noindent{\textbf{Proof. }}}
\def\QED{\hfill {$\square$}\goodbreak \medskip}

\newtheorem{Theorem}{Theorem}[section]
\newtheorem{Lemma}[Theorem]{Lemma}
\newtheorem{Proposition}[Theorem]{Proposition}
\newtheorem{Corollary}[Theorem]{Corollary}
\newtheorem{Remark}[Theorem]{Remark}

\begin{document}

\title 
{Embedded loops in the hyperbolic plane with prescribed, almost constant curvature}
\author{Roberta Musina\footnote{Dipartimento di Scienze Matematiche, Informatiche e Fisiche, Universit\`a di Udine,
via delle Scienze, 206 -- 33100 Udine, Italy. Email: \texttt{roberta.musina{@}uniud.it}. 
{Partially supported by Miur-PRIN project 2015KB9WPT\_001.}
}
\and 
Fabio Zuddas
\footnote{Dipartimento di Matematica e Informatica, Universit\`a di Cagliari, Palazzo delle Scienze, Via Ospedale 72 -- 09124 Cagliari, Italy. Email: \texttt{fabio.zuddas{@}unica.it}.
}}
\date{}

\maketitle

\begin{abstract}
\noindent {Given a constant $k>1$ and
a real valued function $K$  on the hyperbolic plane $\H^2$, we study the problem
of finding, for any $\eps\approx 0$, a closed and embedded curve $u^\eps$ in $\H^2$
having geodesic curvature $k+\eps K(u^\eps)$ at each point.
}
\end{abstract}

\section{Introduction}
Let $\Sigma$ be an oriented Riemannian surface with empty boundary, 
Riemannian metric tensor $g$ and Levi-Civita connection $\nabla^\Sigma$. 
The geodesic curvature of a regular loop $u\in C^2(\S^1,\Sigma)$ is given by
$$
K(u)=\frac{\langle \nabla^\Sigma_{\!u'}u',i_uu'\rangle_g}{|u'|^3_g}\ms.
$$
Here we denoted by $i_u:T_u\Sigma\to T_u\Sigma$  the isometry that rotates 
$T_u\Sigma$,
in such a way that $\{\tau,i_u\tau\}$ is a positively oriented orthogonal basis of $T_u\Sigma$, for any $\tau\neq 0$. 

Given a sufficiently
smooth function $K:\Sigma\to \R$, the {\em $K$-loop problem} consists in
finding 
regular curves $u\in C^2(\S^1,\Sigma)$ having 
geodesic curvature $K(u)$ at each point. This problem can be faced
by studying the system of ordinary differential equations
\begin{equation}
\label{eq:problem_general}
\nabla^\Sigma_{\!u'}u'=L^{\!\Sigma\!}(u)K(u)\ms i_uu'\ms,\quad L^{\!\Sigma\!}(u):=\big(\fintm|u'|^2_g\ms dx\big)^\frac12\ms.
\end{equation}
Indeed, every nonconstant
solution $u\in C^2(\S^1,\Sigma)$ to (\ref{eq:problem_general})
has constant speed $|u'|_g=L^{\!\Sigma\!}(u)$, use for instance the computations in  \cite[Chapter 4]{Jos}. Therefore
$u$ is regular, and has curvature
$K(u)$ at each point.  

The $K$-loop problem has been largely studied since the seminal work \cite{A1} by 
Arnol'd. Most of the
available existence results require compact target surfaces $\Sigma$; we limit ourselves to cite
 \cite{Co, G, G2,NT, S0, S, Sl, T2} and references therein.

\medskip

In the present paper we take $\Sigma$ to be the (noncompact) hyperbolic plane 
$\H^2$. 
It turns out that the problem under consideration does not have solutions,
in general (see Subsection \ref{S:none}). In particular, if 
$-1\le K(q)\le 1$ for any $q\in\Sigma$, then no $K$-loop exists.
If $K\equiv k>1$ is constant (recall that changing the orientation of a  curve  
changes the sign of its curvature), then any regular parameterization of an 
 hyperbolic circle of radius 
\begin{gather*}
\rho_k=\artanh \frac1k =
\frac12\ln\frac{k+1}{k-1}
\end{gather*}
is  a $k$-loop; conversely, 
any  $k$-loop in $\H^2$ parameterizes some circle of radius $\rho_k$.

\medskip

Our existence results involve curvatures that are small perturbations of
a given constant $k>1$. 
In Section \ref{S:constant} we carefully choose 
a reference parameterization $\omega$
of a circle of radius $\rho_k$. Then we take 
any point $z\in\H^2$ and compose $\omega$
with an hyperbolic translation to obtain
a parameterization 
$\omega_z$ of $\partial D^\H_{\rho_k}(z)$. 
Next, given $K\in C^1(\H^2)$, we look for 
a point $z_0\in\H^2$ and for embedded $(k+\eps K)$-loops in $\H^2$ that  suitably approach the circle 
$\omega_{z_0}$ as $\eps \to 0$.

The center $z_0$ can not be arbitrarily prescribed. 
In fact, in Theorem \ref{T:necessary}  we prove that if there exists a sequence of 
$(k+\eps_h K)$-loops $u_h$ such that $\eps_h\to 0$ and $u_h\to \omega_{z_0}$ suitably,
then $z_0$ is a critical point for the Melnikov-type function
\begin{equation}
\label{eq:Mel}
F^{\!K}_k(z)=\int\limits_{D^\H_{\rho_k}(z)} K(z)dV_{\H}~,\qquad F^{\!K}_k:\H^2\to\R~\!.
\end{equation}

One may wonder whether the existence of a critical point $z_0$ for $F^{\!K}_k$ is sufficient to have the existence,
for $\eps\approx 0$, of an embedded $(k+\eps K)$-loop $u_\eps\approx \omega_{z_0}$.
We can give a positive answer in case $F^{\!K}_k$ has a {\em stable} critical point,
accordingly with the next definition (see also \cite[Chapter 2]{AM}).

\bigskip
\noindent
{\bf Definition}
{\em Let ${X}\in C^1(\H^2)$ and let $A\Subset \H^2$ be an open set. We say that ${X}$ has a stable critical point in $A$ if
there exists $r>0$ such that any  function ${G}\in C^1(\overline A)$
satisfying 
$\displaystyle{\|{G}-{X}\|_{C^1(\overline A)}<r}$  has a critical point in $A$.}

\bigskip

Sufficient conditions to have the existence of a stable critical point $z\in A$ for ${X}$
are easily given via elementary calculus. 
For instance, one can assume that one of the following conditions holds:
\begin{itemize}
\item[$i)$] $\nabla {X}(z)\neq 0$ for any $z\in\partial A$, and $\deg(\nabla {X},A, 0)\neq 0$, 
where ''$\deg$'' is Browder's topological degree;
\item[$ii)$] $\displaystyle{{\min_{\partial A} {X}>\min_{A} {X}}}$ ~or~
$\displaystyle{\max_{\partial A} {X}<\max_{A} {X}}$;
\item[$iii)$] ${X}$ is of class $C^2$ on $A$, it has a critical point $z_0\in A$, 
and the Hessian matrix of ${X}$ at $z_0$ is invertible.
\end{itemize}

\medskip

We are in position to state our main  result.

\begin{Theorem}
\label{T:main}
Let $k>1$ and $K\in C^1(\H^2)$ be given. Assume  that $F^{\!K}_k$ has a stable 
critical point in an open set $A\Subset \H^2$. 
Then for every $\eps\in\R$ close enough to $0$, there exists an embedded $(k+\eps K)$-loop $u^\eps$.

Moreover, any sequence $\eps_h\to 0$ has a subsequence $\eps_{h_j}$ such that 
$u^{\eps_{h_j}}\to \omega_{z_0}$ in $C^2(\S^1,\H^2)$ as $j\to\infty$, where
{$z_0\in \overline A$} is a critical point for $F^{\!K}_k$.
In particular, if 
a point $z_0\in A$ is the unique critical point for $F^{\!K}_k$ in $\overline A$, then  
$u^\eps\to \omega_{z_0}$ 
in $C^2(\S^1,\H^2)$ as $\eps\to 0$.
\end{Theorem} 

Any stable critical point of the perturbation term $K$ gives rise to a stable critical point
for $F^{\!K}_k$, at least for $k$ large enough. 
This is in essence the argument we use in Theorem \ref{T:main2}
to obtain, via Theorem \ref{T:main}, the existence of $k+\eps K$-loops 
whenever the perturbation curvature $K$ admits stable critical points.

The proof of Theorem \ref{T:main}
is based on a Lyapunov-Schmidt reduction technique combined 
with variational arguments, as proposed in \cite{AB} (see also \cite[Chapter 2]{AM}).

In fact, $(k+\eps K)$-loops correspond to  critical points of an
energy functional $E_{k+\eps K}(u)=E_{k+\eps K}(u)$, where $u$ runs in the class of nonconstant curves in $C^2(\S^1,\H^2)$
(see Section \ref{S:variational} for details).
In particular, critical points of the unperturbed functional $E_k$ are  circles of
radius $\rho_k$. Let $\mathcal S=\{\omega_z\circ \xi\}$,
where $\xi$ is a rotation of $\S^1$,
$z\in\H^2$, and
$\omega_z$ is our reference parameterization of $\partial D^\H_{\rho_k}(z)$. 
Clearly 
$\mathcal S$ is a smooth three-dimensional manifold of solutions to the unperturbed problem 
$E'_k(u)=0$.

The crucial and technically difficult nondegeneracy result is proved in Lemma \ref{L:nondegen}, via an efficient 
functional change inspired by \cite{Mu}. 
It states that for any $z\in\H^2$, the tangent space to 
$\mathcal S$ at $\omega_z$ coincides with the set of solutions to the linear problem  
$E''_k(\omega_z)\f=0$. 
In the last section we carry out the dimensional reduction argument and complete the proof of Theorem
\ref{T:main}. 

We conclude the paper with a short appendix about the much more easy problem of finding
loops in $\R^2$ having prescribed, almost constant curvature.

The Lyapunov-Schmidt reduction argument has been successfully used to study 
related geometrical problems. We limit ourselves to cite the pioneering paper
\cite{Ye} by R. Ye,
\cite{ALM, Cal, Duke, CM, Fa, Fe, Mo, Mu} and references therein.

\section{Notation and preliminaries}
\label{S:preliminaries}

The Euclidean space $\R^2$ is endowed with the scalar product
$p\cdot q$ and norm $|\cdot|$, so that the disk of radius $R$ centered at $p\in\R^2$ is
$D_R({p})=\{z\in\R^2~|~|z-p|<R\}$. 
The canonical basis of $\R^2$ is $e_1=(1,0), e_2=(0,1)$. 

Let $A,\Omega\subseteq\R^2$ be open sets. We write $A\Subset\Omega$ if 
$\overline A$ is a compact subset of $\Omega$.

\smallskip

We will often use complex notation for points in $\R^2$. In particular
we write $iz=(-z_2,z_1)$ and  $z^2=(z_1^2-z_2^2,2z_1z_2)$ for $z=(z_1,z_2)\in \R^2$.

Let $\S^1$  be the unit circle in the complex plane. Any $\xi\in\S^1$ is identified
with the rotation $x\mapsto \xi x$.

\paragraph{The Poincar\'e half-plane model}~\\
We adopt
as model for the  two dimensional hyperbolic space the half-plane
$$\H^2=\{ z=(z_1,z_2)\in\R^2\ | \ z_2 > 0 \}$$
endowed with the Riemannian metric $g_{lj}(z)= z_2^{-2}\delta_{lj}$.  
With some abuse of notation, we   use the symbol $\H^2$  to denote
the Euclidean upper half space as well.

The hyperbolic distance $d_\H(p,q)$ in $\H^2$ is related to the Euclidean one by 
$$\cosh d_\H(p,q) = 1 + \frac{|p-q|^2}{2 p_2 q_2}\ms,$$ 
and the hyperbolic disk $D^\H_\rho({p})$ centered at $p=(p_1, p_2)$ 
is the Euclidean disk of center 
$(p_1, p_2 \cosh \rho)$ and radius $p_2 \sinh \rho$.

\medskip

A {\em loop} in the $2$-dimensional hyperbolic space $\H^2$ is a
curve $u:\S^1\to\H^2$ of class $C^2$ having nonzero derivative
at each point. We say that $u$ is embedded if it is injective.

\medskip

If $G:\H^2\to\R$ is a differentiable function, then 
$\nabla^{\!\H} G(z)=z_2^2\nabla G(z)$, where $\nabla^{\!\H}$, $\nabla$
are the hyperbolic and the Euclidean gradients, respectively. In particular, 
$\nabla^{\!\H} G(z)=0$ if and only if $\nabla G(z)=0$.

The hyperbolic
volume form $dV_{\H}$ is related to the Euclidean one by $dV_{\H}=z_2^{-2}dz$.

The Levi-Civita connection in $\H^2$ along a curve $u$ in $\Sigma$ is given by
\begin{equation}
\label{eq:nablaHyp}
\nabla^\H_{\!u'}u'= u''-u_2^{-1}\Gamma(u')\ms,
\end{equation}
where, in complex notation, $\Gamma(z)=-iz^2$. In coordinates we have
\begin{equation}
\Gamma(z):=(2z_1z_2,z_2^2-z_1^2)=  z_2z-z_1~\!iz~,\qquad \Gamma:  \H^2\to\R^2~\!.
\label{eq:useful0}
\end{equation}
For future convenience we compute the differential
\begin{equation}
\Gamma'(z)w=2(w_2z-w_1~\!iz)~,\qquad z\in\H^2,~w\in\R^2. 
\label{eq:usefulX}
\end{equation}

\paragraph{Isometries in $\H^2$}~\\
{\em Hyperbolic translations} are obtained by composing a horizontal (Euclidean) translation
$w\mapsto w+se_1$, $s\in\R$ (sometimes called {\em parabolic isometry}), 
with an Euclidean homothety $w\mapsto tw$, $t>0$ (in some literature,
only homotheties are called  hyperbolic translations). 
We obtain the two dimensional group of isometries  $\H^2\to\H^2$,
$$
u\mapsto u_z:=z_1e_1+z_2u~,\qquad z\in\H^2\ms.
$$

\paragraph{Function spaces}~\\
Let $m\ge 0$, $n\ge 1$ be integer numbers. We endow $C^m(\S^1,\R^n)$ with the standard Banach space structure.
If $f\in C^1(\S^1,\R^n)$, we identify $f'(x)\equiv f'(x)(ix)$, so that $f':\S^1\to\R^n$.

In  $L^2=L^2(\S^1,\R^2)$ we take the Hilbertian norm 
$$\displaystyle{\|u\|^2_{L^2}=\fintm |u(x)|^2~\!dx~\!=\frac{1}{2\pi} \irn|u(x)|^2~\!dx}\ms.$$
If $T\subseteq C^0(\S^1,\R^2)$, the orthogonal to $T$ with respect to the $L^2$ scalar product is 
$$
T^\perp=\{\f\in C^0(\S^1,\R^2)~|~\fintm u\cdot\f~\!dx=0~~\text{for any $u\in T$}~\}.
$$

We look at $C^m(\S^1,\H^2)$ as an open subset
of the Banach space $C^m(\S^1,\R^2)$, and  identify $\H^2$ with the set of constant functions in 
$C^m(\S^1,\H^2)$. Thus 
$C^m(\S^1,\H^2)\setminus\H^2$ contains only nonconstant curves.

\subsection{The  variational approach}
\label{S:variational}

We put
$$
L(u):=L_{\H^2}(u)=\Big(\fintm u_2^{-2}|u'|^2~\!dx\Big)^\frac{1}{2}~,\qquad L:C^2(\S^1,\H^2)\to\R\ms,
$$
that is a $C^\infty$ functional, with  Fr\'echet differential
\begin{equation}
\label{eq:dL0}
L'(u)\f
=\frac{1}{L(u)}~\!\fintm {u_2^{-2}}\big(-u''+u_2^{-1}\Gamma(u')\big)\cdot\f~\!dx\ms,
\quad \f\in C^2(\S^1,\R^2)\ms.
\end{equation}
When $\Sigma=\H^2$, problem (\ref{eq:problem_general})   reads
\begin{equation*}\tag{$\mathcal P_K$}
\label{eq:EL}
u''-u_2^{-1}\Gamma(u')=L(u)K(u)~\!iu'\ms.
\end{equation*}
The system (\ref{eq:EL}) admits a variational formulation.  
More precisely,
its nonconstant
solutions are  critical points of the energy functional of the form
$$
E_K(u)= L(u)+A_K(u)~,\quad u\in C^2(\S^1,\H^2)\setminus\H^2,
$$
where $A_K(u)$  gives, roughly speaking,
the signed area enclosed by the curve $u$ with respect to the weight $K$
(see  Remark \ref{R:area} below). More precisely, to introduce $A_K(u)$ we take any vectorfield
$Q_K\in C^1(\H^2,\R^2)$ such that 
$$
\div Q_K(z)=z_2^{-2}K(z)~,\qquad z\in\H^2
$$
(here $''\div''$ is the usual Euclidean divergence). A possible choice is
$$
Q_K(z_1,z_2)=\Big(\ms\frac12\ms z_2^{-2}\int\limits_0^{z_1}K(t,z_2)\ms dt\Big)e_1+
\Big(\ms\frac12\ms \int\limits_1^{z_2}t^{-2}K(z_1,t)\ms dt\Big)e_2\ms.
$$
Then we define
$$
A_K(u)=\fintm Q_K(u)\cdot iu'~\!dx~,\qquad A_K:C^2(\S^1,\H^2)\to\R.
$$
By direct computations one gets that the functional $A_K$ is 
Fr\'echet differentiable at any $u\in C^2(\S^1,\H^2)$, with
differential  
\begin{equation}
\label{eq:dAK}
A_K'(u)\f=~\!\fintm {u_2^{-2}}\ms K(u)\f \cdot\ms iu'~\!dx\ms.
\quad \f\in C^2(\S^1,\R^2)\ms,
\end{equation}
It follows that 
$A_K(u)$ does not depend on the choice of the vectorfield $Q_K$. Further,
if $K\in C^1(\H^2)$ then 
the area functional $A_K$ is of class $C^2$ on $C^2(\S^1,\R^2)$.

In conclusion, the following lemma holds. 
\begin{Lemma}
\label{L:E}
Let $K\in C^1(\H^2)$. The functional 
$E_K(u)=L(u)+A_K(u)$
is of class $C^2$ on $C^2(\S^1,\H^2)\setminus\H^2$, and
$$
L(u)E'_K(u)\f= \fintm {u_2^{-2}}\big(-u''+u_2^{-1}\Gamma(u')+L(u)K(u)\ms iu'\big)\cdot\f~\!dx\ms
$$
for any $u\in C^2(\S^1,\H^2)\setminus\H^2, \f \in C^2(\S^1,\R^2)$.
In particular, if $u_0\in C^2(\S^1,\H^2)\setminus\H^2$ is a critical point for the functional
$E_K(u)$, then $u_0$ solves {\rm (\ref{eq:EL})}, hence it is an hyperbolic $K$-loop.
\end{Lemma}

\begin{Remark}
\label{R:area}
Let $u\in C^2(\S^1,\H^2)$ be an embedded loop. 
Then $u$ is a regular parameterization of the boundary of an open set $\Omega_u\Subset\H^2$.
Assume for instance that $u$ is positively oriented, so that
$iu'$ gives the inner direction to $\Omega_u$. Then
$$
A_K(u)=
-\frac{1}{2\pi}\int\limits_{\partial \Omega} Q_K(z)\cdot \nu \ms ds= 
- \frac{1}{2\pi}\int\limits_{\Omega} K(z) dV_\H
$$
by the divergence theorem. 
\end{Remark}

\subsection{Nonexistence results}
\label{S:none}

We start with a simple result that  should be well known. 
We sketch its proof by adapting the argument in \cite[p. 194]{Lop}.

\begin{Proposition}
\label{teoLopez}
Let $K\in C^0(\H^2)$. If $\|K\|_\infty\le 1$ then no $K$-loop exists.
\end{Proposition}

\proof
Let $u\in C^2(\S^1,\H^2)$ be a $K$-loop. We need to show that $|K|>1$ somewhere in $\H^2$.
Take the smallest closed disk $D_\rho=\overline{D^\H_\rho(z)}$ containing $u(\S^1)$. Then
$\partial D_\rho$ is tangent to $u(\S^1)$ at some point. 
At the contact point the absolute value of the curvature of $u$ can not be smaller than the curvature 
$1/\tanh\rho$ of the circle $\partial D_\rho$, use a local comparison principle. 
The conclusion readily follows from $\tanh\rho<1$.
\QED

Next, we point out few necessary conditions for the existence of $K$-loops.

\begin{Lemma}
\label{L:none}
Let $K\in C^1(\H^2)$ and let
$\Omega \subset \H^2$ be a bounded open domain. Assume that $\partial\Omega$
is parameterized by a $K$-loop $u\in C^2(\S^1,\H^2)$.
Then
$$
\int\limits_{\Omega} \nabla K(z)\cdot e_1\ms   dV_{\H^2}=
\int\limits_{\Omega} \nabla K(z)\cdot z\ms   dV_{\H^2}=
\int\limits_{\Omega} \nabla K(z)\cdot z^2\ms   dV_{\H^2}=0\ms.
$$
\end{Lemma}

\proof
Direct computations based on integration by parts give
\begin{equation}
\label{eq:kerL}
L'(u)e_1=L'(u)u=L'(u)i(\Gamma u)=0,
\end{equation}
see (\ref{eq:dL0}) and (\ref{eq:useful0}).
In addition, the curve $u$ solves
$$
-L(u)L'(u)\f=\fintm u_2^{-2}K(u)\f\cdot iu'\ms dx\quad\text{for any $\f\in C^2(\S^1,\R^2)$}.
$$
Since $iu'(x)\neq 0$ is parallel to the outer normal $\nu$ to $\Omega$ at $u(x)\in\partial\Omega$, we infer that
$$
\int\limits_{\partial\Omega} z_2^{-2}K(z)e_1\cdot \nu\ms =\int\limits_{\partial\Omega} z_2^{-2}
K(z)z\cdot \nu\ms =\int\limits_{\partial\Omega} z_2^{-2} K(z)i\Gamma(z)\cdot \nu\ms =0.
$$
Recall that we identify $i\Gamma(z)=z^2$, then use the divergence theorem to get
$$
\int\limits_{\Omega} \div\big(z_2^{-2}K(z)e_1\big)\ms dz=\int\limits_{\Omega}  \div\big(z_2^{-2}K(z)z\big)\ms dz
=\int\limits_{\Omega}  \div\big(z_2^{-2}K(z)z^2\big)\ms dz=0\ms.
$$
The conclusion readily follows.
\QED

\begin{Remark}
The  identities in (\ref{eq:kerL})  hold indeed for any curve $u$, and are related to the group of isometries in $\H^2$.
Notice indeed that $z\mapsto e_1, z\mapsto z, z\mapsto z^2$ are infinitesimal Killing vectorfields in $\H^2$.
\end{Remark}

Lemma \ref{L:none} readily implies the next nonexistence result.

\begin{Corollary}
\label{C:none}
Let $K\in C^1(\H^2)$ be a given curvature function. 
Assume that one of the following conditions hold,
\begin{itemize}
\item[$i)$] $K$ is  
strictly monotone in the $e_1$ direction;
\item[$ii)$] $K$ is radially 
strictly monotone, that is, $\nabla K(z)\cdot z$ never vanishes on $\H^2$;
\item[$iii)$] $\nabla K(z)\cdot z^2$ never vanishes on $\H^2$
\end{itemize}
Then no embedded $K$-loop exists.
\end{Corollary}

\section{The unperturbed problem}
\label{S:constant}

\small 
In this section we take a constant $k>1$ and study the system 
\begin{equation*}\tag{$\mathcal P_k$}
\label{eq:kequ}
u''-u_2^{-1}\Gamma(u')={L(u)}k~\!iu'\ms. 
\end{equation*}
We start by introducing the radius
$$
{R_k}:=\sinh\rho_k=\frac{1}{k}\cosh \rho_k=\frac{1}{\sqrt{k^2-1}}
$$
and the reference loop $\omega:\S^1\to\H^2$,
\begin{equation}
\label{eq:def_omega}
\omega(x)=\frac{1}{k-x_2}\big(\ms x_1\ms ,\frac{1}{R_k}~\!\big)~,\quad x=x_1+ix_2\in\S^1\ms.
\end{equation}
Notice that 
\begin{equation}
\label{eq:omega1}
|\omega-kR_k e_2|={R_k},
\end{equation}
hence $\omega$ is a (positive) parametrization of the Euclidean circle $\partial D_{R_k}(kR_k e_2)$, that 
coincides with the
hyperbolic circle $\partial D_{\rho_k}^\H(e_2)$.  
The next identities will be very useful:
\begin{eqnarray}
\label{eq:omega2}
&&\omega'=\omega_2~\!i(\omega-kR_k e_2)\\
&&
\label{eq:3.4}
\omega_2^{-1}\Gamma(\omega')=(\omega_2-kR_k)~\!i\omega'+\omega_1~\!\omega'\\
\label{eq:omega3}
&&\omega_2^{-1}|\omega'|\equiv L(\omega)={R_k}\ms .
\label{eq:omega4}
\end{eqnarray}
By differentiating (\ref{eq:omega2}) and using (\ref{eq:omega4}) one easily gets that $\omega$ 
solves (\ref{eq:kequ}). Next, for $z=(z_1,z_2)\in\H^2$ we  parameterize $\partial D^\H_{\rho_k}(z)$ by the function
$$
\omega_z=z_1e_1+z_2\omega~\!.
$$
Notice that $\omega=\omega_{e_2}$. It is easy to check that for any 
rotation $\xi\in\S^1$ and any point $z\in\H^2$,
the circle $\omega_z\circ \xi$  solves (\ref{eq:kequ}) as well. 
Further, by Remark \ref{R:area} we have
\begin{equation}
\label{eq:pioggia}
F^{\!K}_k(z):=\int\limits_{D^\H_{\rho_k}(z)} K(z)dV_{\H}=-2\pi A_K(\omega_z).
\end{equation}

We know that any nonconstant solution $u$ to (\ref{eq:kequ}) has
constant curvature $k$, hence is a circle of hyperbolic radius $\rho_k$. Actually we
need a sharper uniqueness result, that is,
we have to classify solutions to (\ref{eq:kequ}).

\begin{Lemma}
\label{L:k_uniqueness}
Let $u\in C^2(\S^1,\H^2)$ be a nonconstant solution to (\ref{eq:kequ}). Then
$\mu:=\displaystyle{L(u)}/{L(\omega)}$ is an integer number, and
 there exist $\xi\in\S^1$, $z=(z_1,z_2)\in\H^2$  such that 
$u(x)= \omega_z\circ \xi$. 
In particular, $u$ parameterizes $\partial D_{\rho_k}(z)$, and 
$ L(u)=\mu L(\omega)=\mu R_k$.
\end{Lemma}

\proof
We have
$$
\omega_2(-i)=e^{-\rho_k}=\min_{x\in\S^1}\omega_2(x)\ms~,\quad
\omega(-i)= e^{-\rho_k} e_2~,\quad \omega'(-i)= e^{-\rho_k} L(\omega) e_1\ms.
$$
Let $x_u\in\S^1$ such that 
$$
u_2(x_u)=m_u:=\min_{x\in\S^1} u_2(x)\ms.
$$
Now we show that
\begin{equation}
\label{eq:claimm}
u'(x_u)= m_uL(u) e_1.
\end{equation}
Clearly $u'_2(x_u)=0$ and $u''_2(x_u)\ge 0$. We first infer that 
$\Gamma(u'(x_u))=-u'_1(x_u)~\!iu'(x_u)$, compare with (\ref{eq:useful0}). Thus
the system (\ref{eq:kequ}) for the second coordinate
gives
$$
\big(L(u)k-m_u^{-1}u'_1(x_u)\big)~\!u'_1(x_u)=u_2''(x_u)\ge 0,
$$
that implies $u'_1(x_u)\ge 0$. On the other hand, 
$u_2^{-1} |u'|\equiv L(u)$ on $\S^1$. 
Thus $u'_1(x_u)=|u'(x_u)|=m_u L(u)$, and (\ref{eq:claimm}) is proved.

In particular, 
$u$ solves the Cauchy problem
\begin{equation}
\label{eq:cauchy}
v''= v_2^{-1}\Gamma(v')+kL(u)~\!i v'~,\quad
v(x_u)=u(x_u)~,\quad
v'(x_u)= m_uL(u) e_1.
\end{equation}
It is easy to check that the function
$$
\tilde u(x):={m_u}{e^{\rho_k}}~
\!\omega\big(-ix_u^{-\mu}x^\mu\big)+u_1(x_u)e_1
$$
solves (\ref{eq:cauchy}) as well (use $f'(x)=i\mu x^\mu$ for $f(x)=x^\mu$, $f:\S^1\to\C$). 
Thus $\tilde u(x)= u(x)$ for any $x\in\S^1$
and hence $u(x)= \omega_z\circ \xi$, where 
$z_1=u_1(x_u)$, $z_2={m_u}{e^{\rho_k}}$, $\xi= -ix_u^{-\mu}$.
Finally, $\mu$ is an integer number because $u$ and $\omega$ are both well defined on $\S^1$.
\QED

\paragraph{The linearized  problem}~\\
By Lemma \ref{L:k_uniqueness},  the $3$-dimensional manifold
$$
\mathcal S=\big\{\omega_z\circ \xi~|~\xi\in\S^1~,~~z\in\H^2~\big\}\subset C^2(\S^1,\H^2),
\quad \omega_z= z_1e_1+z_2\omega
$$
is  the set of embedded solutions to (\ref{eq:kequ}). 
The tangent space to $\mathcal S$ at $\omega_z$ is
$$
T_{\omega_z}\mathcal S=T_{\omega}\mathcal S=\langle ~\omega',e_1,\omega~\rangle.
$$
Every loop in $\omega_z\circ \xi\in \mathcal S$ is a critical point for the energy functional 
$$
E_k(u)=L(u)+A_k(u)=\Big(\fintm u_2^{-2}|u'|^2~\!dx\Big)^\frac{1}{2}-k\fintm u_2^{-1} u'_1\ms dx
$$
on $C^2(\S^1,\H^2)\setminus H^2$, and $E_k(\omega_z\circ \xi)=E_k(\omega)$ is a constant. More generally one has
\begin{equation}
\label{eq:invariance0}
E_k(z_1e_1+z_2\ms u\circ\xi)=E_k(u)\qquad\text{for any $ \xi\in\S^1$, $z\in\H^2$.}
\end{equation}

In order to handle the differential of $E_k$, it is convenient to introduce the function $J_0:C^2(\S^1,\H^2)\setminus \H^2\to C^0(\S^1,\H^2)$ given by
\begin{eqnarray}
\nonumber
J_0(u)&=& -(u_2^{-2}u')'-u_2^{-3}|u'|^2e_2+L(u)k u_2^{-2} ~\!iu'\\
&=&
u_2^{-2}\big(-u''+u_2^{-1}\Gamma(u')+L(u)k\ms iu'\big)\ms.
\label{eq:J0}
\end{eqnarray}
By Lemma \ref{L:E} we have
\begin{equation}
\label{eq:dL}
L(u)E'_k(u)\f 
=\fintm J_0(u)\cdot\f~\!dx\quad \text{for any $\f \in C^2(\S^1,\R^2)$.}
\end{equation}
By differentiating  (\ref{eq:invariance0}) at $\xi=1$, $z=e_2$ we readily get
$E'_k(u)u'=E'_k(u) e_1=E'_k(u)u=0$ for any nonconstant curve 
$u\in C^2(\S^1,\H^2)$, that is, 
\begin{equation}
\label{eq:inv0}
\fintm J_0(u)\cdot u'~\!dx=0~,\quad \fintm J_0(u)\cdot e_1~\!dx=0~,\quad
\fintm J_0(u)\cdot u~\!dx=0\ms.
\end{equation}
Now we differentiate (\ref{eq:dL}) with respect to $u$, at $u=\omega_z$. From $E'_k(\omega_z)=0$ we get
$$
L(\omega)E''_k(\omega_z)[\f,\tilde\f]=\fintm J_0'(\omega_z)\f\cdot\tilde\f~\!dx \quad \text{for any $\f,\tilde\f\in C^2(\S^1,\R^2)$.}
$$
Since $E_k$ is of class $C^2$, then $J'_0(\omega_z)$ is self-adjoint in $L^2$, that means
\begin{equation}
\label{eq:autoaggiunto}
\fintm J_0'(\omega_z)\f\cdot\tilde\f~\!dx=\fintm J_0'(\omega_z)\tilde\f\cdot\f~\!dx
\quad \text{for any $\f,\tilde\f\in C^2(\S^1,\R^2)$.}
\end{equation}
Finally, we differentiate $E'_k(\omega_z\circ\xi)=0$
with respect to the variables $\xi\in\S^1,
z\in\H^2$ to get  $T_{\omega_z}\mathcal S\subseteq \ker\! J_0'(\omega_z)$.
We shall see in the crucial Lemma \ref{L:nondegen} below that indeed 
$T_{\omega_z}\mathcal S= \ker\! J_0'(\omega_z)$. 

This will be done via a useful functional change.

\paragraph{A functional change and nondegeneracy}~\\
In order to avoid tricky computations, we use 
in $C^m(\S^1,\R^2)$, $m\ge 0$, the orthogonal frame $\omega', i\omega'$.
We introduce the isomorphism 
$$
\Phi(g)=g_1 \omega'+g_2 i\omega'~,\qquad \Phi:C^m(\S^1,\R^2)\to C^m(\S^1,\R^2)
$$
together with its inverse $\Phi^{-1}(\f)=R_k^{-2}\omega_2^{-2}(\ms \f\cdot \omega'\ms e_1+\f\cdot i\omega'\ms e_2)$ (recall
that $|\omega'|=R_k\omega_2$) and
the differential operator 
\begin{equation}
\label{eq:B}
Bg=
-g''-{kR_k} ig'+R_k^2 \big(g_2-k^2\fintm g_2dx\big)e_2~,\qquad g\in C^2(\S^1,\R^2)\ms.
\end{equation}

\begin{Lemma}
\label{L:new}
Let $z$ be any point in $\H^2$. The following facts hold.
\begin{itemize}
\item[$i)$] $J_0'(\omega_z)(\Phi(g))=z_2^{-2}\omega_2^{-2}\Phi\big(Bg\big)$
for any $g\in C^2(\S^1,\S^2)$;
\item[$ii)$] 
$\displaystyle{\fintm\omega_2^{-2}\Phi(g)\cdot\Phi(\tilde g)~\!dx=R_k^2\fintm g\cdot \tilde g~\!dx}$
for any $g,\tilde g\in C^2(\S^1,\S^2)$
\end{itemize}
\end{Lemma}

\proof
Since $J_0(\omega_z)=z_2^{-1}J_0(\omega)=0$ and $J_0'(\omega_z)=z_2^{-2}J_0'(\omega)$,
it suffices to prove $i)$ for $z=e_2$, that corresponds to $\omega_z=\omega$.
We have to show that
\begin{equation}
\label{eq:app}
\mathcal J(\f):=\omega_2^2J_0'(\omega)\f=\Phi\big(Bg)~,\quad\text{where}~~ \f=g_1\omega'+g_2~\!i\omega'~\!.
\end{equation}
To compute $\mathcal J(\f)$ it is convenient to recall (\ref{eq:J0}) and to differentiate the identity
$$
u_2^2~\! J_0(u)=- u''+u_2^{-1}\Gamma(u')+L(u)k~\!iu'
$$
at $u=\omega$. Since $J_0(\omega)=0$ and $L(\omega)=R_k$,  we get 
$$
\mathcal J(\f)=
-\f''+kR_k~\!i\f'+\omega_2^{-1}\Gamma'(\omega')\f'
-\omega_2^{-2}\f_2\Gamma(\omega')+k\big(L'(\omega)\f\big)i\omega'.
$$
From (\ref{eq:usefulX}) we find
$\Gamma'(\omega')\f'=2\f'_2~\!\omega'-2\f'_1~\!i\omega'$.
Taking also  (\ref{eq:3.4}) into account, we obtain
\begin{eqnarray}
\nonumber
\mathcal J(\f) 
&=&-\f''+{kR_k}~\!i\f'+A_1(\f)~\!\omega'
-\big(A_2(\f)-k~\! L'(\omega)\f~\!\big)~\!i\omega'~\!,
\label{eq:boring0}
\end{eqnarray}
where
$$
A_1(\f)=\big(2\f'_2-\f_2\omega_1\big)\omega_2^{-1}~,
~~ A_2(\f)=\big(2\f'_1+\f_2(\omega_2-{kR_k})\big)\omega_2^{-1}.
$$
To compute the differential $L'(\omega)$ at $\f$ we recall that 
$\omega$ solves (\ref{eq:kequ}). Thus  (\ref{eq:dL0}) gives
$$
L'(\omega)\f=-k \fintm \omega_2^{-2}\f\cdot i\omega'~\!dx~\!.
$$
For the next computations we observe that the loop $\omega$ solves several useful differential systems. 
In particular, from (\ref{eq:kequ}), (\ref{eq:3.4}), (\ref{eq:omega3}) and (\ref{eq:omega2})
it follows that
\begin{equation}
\label{eq:omega5}
\omega''=\omega_1 \omega'+\omega_2~\!i\omega'~,\quad
\omega'''=(\omega_1^2-2\omega_2+kR_k\omega_2)~\!\omega'+3\omega_1\omega_2~\!i\omega'.
\end{equation}

\medskip

Now we take any $\psi\in C^2(\S^1,\R)$ and we look for an explicit formula for 
$\mathcal J(\psi \omega')$. Clearly $L'(\omega)(\psi \omega')=0$,
as $\omega'\cdot i\omega'\equiv 0$. Direct computations based on (\ref{eq:omega5}) give
\begin{eqnarray*}
-(\psi\ms \omega')''+kR_k\ms i(\psi\ms \omega')'&=&
\big(-\psi''-2\omega_1\psi'-(\omega_1^2-2\omega_2^2+2kR_k\omega_2)\psi\big)\ms \omega'\\
&&\qquad+\big((kR_k-2\omega_2)\psi'+(kR_k-2\omega_2)\omega_1\psi\big)\ms i\omega'
\\
~&&\\
A_1(\psi\ms \omega')&=&2\omega_1\psi'-(2\omega_2^2-2kR_k\omega_2-\omega_1^2)\psi\\
A_2(\psi\ms \omega')&=&2(kR_k-\omega_2)\psi'-(kR_k-3\omega_2)\psi,
\end{eqnarray*}
and we find the formula
\begin{equation}
\label{eq:seconda}
\mathcal J(\psi \omega')=-\psi''\omega'-kR_k\psi'\ms i\omega'.
\end{equation}
Now we handle $\mathcal J(\psi ~\!i\omega')$.  
From (\ref{eq:omega4}) we get
$$
k~\! L'(\omega)(\psi\ms i\omega')=-k^2\fintm \omega_2^{-2}|\omega'|^2\psi~\!dx=-k^2R_k^2\fintm\psi~\!dx.
$$
Then we use (\ref{eq:omega1}--\ref{eq:omega4}) and (\ref{eq:omega5}) to compute
\begin{eqnarray*}
-(\psi\ms i\omega')''+kR_k\ms i(\psi\ms i\omega')'&=&
\big((2\omega_2-kR_k)\psi'+(3\omega_2-kR_k)\omega_1\psi\big)\ms \omega'\\
&&\quad+\big(-\psi''-2\omega_1\psi'-(\omega_1^2-2\omega_2^2+2kR_k\omega_2)\psi\big)\ms i\omega'
\\
~&&\\
A_1(\psi\ms i\omega')&=&-2(\omega_2-kR_k)\psi'-(3\omega_2-kR_k)\omega_1\psi\\
A_2(\psi\ms i\omega')&=&-2\omega_1\psi'+(\omega_2^2-k^2R_k^2-2\omega_1^2)\psi.
\end{eqnarray*}
Since $\displaystyle R_k^2=|\omega-kR_k e_2|^2=|\omega|^2-2kR_k\ms\omega_2+k^2R_k^2$ by (\ref{eq:omega1}),
we arrive at
$$
\mathcal J(\psi\ms i\omega')=kR_k\psi'\ms\omega'+\big(-\psi''+R_k^2\psi-k^2R_k^2\fintm\psi\ms dx\big)\ms i\omega',
$$
that together with (\ref{eq:seconda}) gives 
\begin{eqnarray*}
\mathcal J(g_1\ms \omega'+g_2\ms i\omega')&=&
\big(-g''_1-{kR_k} g'_2\big)\omega'\\
&&\quad+
\big(-g''_2+{kR_k} g'_1+\displaystyle R_k^2g_2-k^2R_k^2 \fintm g_2~\!dx\big)i\omega'
\end{eqnarray*}
and concludes the proof of (\ref{eq:app}). The proof of $i)$ is complete;
the formula in $ii)$ is immediate, because $\omega'\cdot i\omega'\equiv 0$ and $|\omega'|=R_k\omega_2$.
\QED

We are in position to prove the main result of this section.

\begin{Lemma}[Nondegeneracy]
\label{L:nondegen}
Let $z$ be any point in $\H^2$. The following facts hold.
\\
$~~i)$ $\displaystyle{\ker\! J'(\omega_z)=T_{\omega}\mathcal S}$;
\\
$~ii)$ If $J_0'(\omega_z)\f\in T_\omega\mathcal S$, then $\f\in T_\omega\mathcal S$;
\\
$iii)$
For any $u\in T_{\omega}\mathcal S^\perp$
there exists a unique $\f\in C^2(\S^1,\R^2)\cap T_{\omega}\mathcal S^\perp$ such that $J_0'(\omega_z)\f=u$. 
\end{Lemma}

\proof 
We start by studying the kernel of the operator $B$ in (\ref{eq:B}). In coordinates,
the linear problem $Bg=0$ becomes
$$
-g''_1+kR_k g'_2=0~,~~
-g''_2-{kR_k} g'_1+R_k^2 \Big(g_2-k^2\displaystyle\fintm g_2dx\Big)=0\ms,
$$
that is clearly equivalent to
\begin{equation}
\label{eq:Bsys2}
~\displaystyle\fint_{\S^1} g_2dx=0,~~
-g''_1+kR_k g'_2=0~,~~
-g''_2-{kR_k} g'_1+R_k^2 g_2=0
\end{equation}
because $k>1$. The system (\ref{eq:Bsys2}) can be studied  via
elementary techniques. The conclusion is that 
$\ker\! B=\langle e_1,\gamma,\gamma'\rangle$, where 
$\gamma=\frac{1}{R_k} (\ms kR_k x_1\ms, -x_2\big)$.
Since $\Phi(e_1)=\omega'$, $\Phi(\gamma)=\omega$ and $\Phi(\gamma')=e_1-\omega'$, thanks to
Lemma \ref{L:new} we have
$$
\ker\! J_0'(\omega_z)=\Phi(\ker\ms B)=\Phi\big(\langle e_1, \gamma,\gamma'\rangle\big)=T_\omega\mathcal S,
$$
and the first claim is proved.

Now we prove $ii)$. If $\tau:=J_0'(\omega_z)\f\in T_\omega\mathcal S=\ker J'(\omega_z)$,
then $J_0'(\omega_z)\tau=0$. Taking  (\ref{eq:autoaggiunto}) into account, we obtain
$$
\fintm|J_0'(\omega_z)\f|^2~\!dx=\fintm J_0'(\omega_z)\f\cdot\tau~\!dx=\fintm J_0'(\omega_z)\tau\cdot\f~\!dx=0.
$$
Thus $J_0'(\omega_z)\f=0$, that means $\f \in T_\omega\mathcal S$. 

\medskip

It remains to prove $iii)$. If $u\in T_\omega\mathcal S^\perp$, then $\Phi^{-1}(\omega_2^2 u)$ is orthogonal to
$\ker\! B$ by $ii)$ in Lemma \ref{L:new}. One can compute the Fourier coefficients
of the unique solution $g_u\in \ker B^\perp$ of the system $Bg_u=\Phi^{-1}(\omega_2^2 u)$. 
Then $J_0'(\omega)(z_2^2\Phi(g_u))=u$ by $i)$ in Lemma \ref{L:new}. The function 
$\f$ defined as the $L^2$-projection of $z_2^2\Phi(g_u)$ on $ T_\omega\mathcal S^\perp$ solves $J_0'(\omega)\f=u$
as well, and is uniquely determined by $u$.

The lemma is completely proved.
\QED

\section{The perturbed problem}
\label{S:nonconstant}

Let $k>1$, $K\in C^1(\H^2)$ be given, and let $\eps\in\R$ be a varying parameter. In this section we
study the system
\begin{equation*}\tag{$\mathcal P_{k+\eps K}$}
\label{eq:ELeps}
u''-u_2^{-1}\Gamma(u')=L(u)(k+\eps K(u))~\!iu'\ms.
\end{equation*}

We start  with a necessary condition for the existence of solutions to 
(\ref{eq:ELeps}) having some prescribed behavior as $\eps\to 0$.

\begin{Theorem}
\label{T:necessary}
Let $k>1$, $K\in C^1(\H^2)$, and  $\eps_h\to 0$ be given. For any integer 
$h$, let $u_h\in C^2(\S^1,\H^2)\setminus\H^2$ be a solution to 
\begin{equation*}\tag{$\mathcal P_{\eps_h}$}
\label{eq:uh}
u_h''=(u_h)_2^{-1}\Gamma(u'_h)+L(u_h)(k+\eps_hK(u_h))\ms iu_h',
\end{equation*}
and assume that

\medskip

\centerline{$ L(u_h)\to L_\infty>0$,
\qquad $u_h\to U$ uniformly, for some $U\in C^0(\S^1,\H^2)$.}

\medskip

\noindent
Then  
there exist  $\mu\in\mathbb N$, $\xi\in\S^1$ and a critical point
$z\in\mathbb H^2$ for $F^{\!K}_k$, such that
$U(x)=\omega_z\big(\xi x^\mu)$.
\end{Theorem}

\proof
We have  $|u'_h|\equiv L(u_h)(u_h)_2$, thus 
the sequence $|u'_h|$ is uniformly bounded.
It follows that $u''_h$ is uniformly bounded as well, because $u_h$ solves
(\ref{eq:uh}). Thus, $u'_h$ is bounded in $C^{0,{s}}$ for any ${s}\in(0,1)$
and using (\ref{eq:uh}) again we infer that the sequence $u_h$ converges in $C^{2,{s}}$
 for any ${s}\in(0,1)$. In particular, $U\in C^{2}(\S^1,\H^2)$, $L_\infty=L(U)$
and $U$ solves
$$
U''=U_2^{-1}\Gamma(U')+L(U)k\ms iU'.
$$
Lemma \ref{L:k_uniqueness} applies and gives the existence of 
$\xi\in\S^1$, $z\in\H^2$, $\mu\in\mathbb N$ such that 
$U(x)=\omega_z(\xi x^\mu)$ and  $L_\infty=L(U)=\mu L(\omega)$.

It remains to prove that $z$ is a critical point for $F^{\!K}_k$. We rewrite (\ref{eq:uh}) in the form
\begin{equation}
\label{eq:uhh}
J_0(u_h)+ \eps_hL(u_h) (u_h)_2^{-2}K(u_h)\ms iu_h'=0,
\end{equation}
see (\ref{eq:J0}). Then we test (\ref{eq:uhh})
with the functions $e_1$ and $u_h$. Taking
(\ref{eq:inv0}) into account, we find
$$
\fintm (u_h)_2^{-2} K(u_h) ~\!e_1\cdot iu_h'~\!dx=0~,\quad\fintm  (u_h)_2^{-2}K(u_h)~\!u_h\cdot iu_h'~\!dx=0.
$$
Since $u_h\to U(x)=\omega_z(\xi x^\mu)$, in the limit as $h\to\infty$  we obtain 
$$
{\mu} \fintm (\omega_z)_2^{-2} K(\omega_z) ~\!e_1\cdot i\omega_z'~\!dx=0~,\quad
{\mu} \fintm (\omega_z)_2^{-2} K(\omega_z) ~\!\omega_z\cdot i\omega_z'~\!dx=0\ms,
$$
that is,
$$
\partial_{z_1}A_K(\omega_z)=A'_K(\omega_z)e_1=0~,\quad
\partial_{z_2}A_K(\omega_z)=A'_K(\omega_z)\omega=0.
$$
Thus $z$ is a critical point for $F^{\!K}_k$ because of (\ref{eq:pioggia}). 
\QED

\subsection{Finite dimensional reduction}
By Lemma
\ref{L:E}, $k+\eps K$-loops are the critical points of the functional
$$
E_{k+\eps K}(u)=E_{k}(u)+\eps A_{K}(u)=L(u)+kA_1(u)+\eps A_{K}(u)~,\quad u\in C^2(\S^1,\H^2)\setminus\H^2~\!.
$$
We introduce the 
$C^1$ function $J_\eps: C^2(\R,\H^2)\setminus\H^2\to C^0(\R,\H^2)$,
\begin{eqnarray*}
J_\eps(u)&=&J_0(u)+\eps L(u)u_2^{-2}K(u)~\!iu'\\
&=& 
u_2^{-2}\big(-u''+u_2^{-1}\Gamma(u')+L(u)(k+\eps K(u))\ms iu'\big)\ms,
\end{eqnarray*}
compare with (\ref{eq:J0}), so that
\begin{equation}
\label{eq:dL_eps}
L(u)\ms E'_{k+\eps K}(u)\f=\fintm J_\eps(u)\cdot\f~\!dx\ms~,\quad u\in C^2(\S^1,\H^2), ~\f\in C^2(\S^1,\R^2).
\end{equation}
We will look for critical points for $E_{k+\eps K}$ by solving the problem $J_\eps(u)=0$.

\medskip

First, we notice that  $E_{k+\eps K}(u\circ\xi)= E_{k+\eps K}(u)$ for any  $\xi\in\S^1$, that implies
\begin{equation}
\label{eq:eps_invariance}
\fintm J_\eps(u)\cdot u'~\!dx=0
\qquad\text{for any $\eps\in \R$, $u\in C^2(\S^1,\H^2)\setminus\H^2$.}
\end{equation}

In the next crucial lemma we carry out the Lyapunov-Schmidt procedure, 
in which we take advantage of the variational structure of problem (\ref{eq:ELeps}).

\begin{Lemma}
\label{L:reduction}
Let $\Omega\Subset \H^2$ be a given open set.
There exist $\beps>0$ and  a $C^1$ function
\begin{gather*}
[-\beps,\beps]\times\overline\Omega\to C^2(\S^1,\H^2)\setminus\H^2~,\quad
(\eps,z)\mapsto u^\eps_z
\end{gather*}
such that the following facts hold.
\begin{itemize}
\item[$i)$] $u^\eps_z$ is an embedded loop and  $u^0_z=\omega_z$;
\item[$ii)$] $u^\eps_z-\omega_z\in T_\omega\mathcal S^\perp$;
\item[$iii)$] $J_\eps(u^\eps_z) \in T_\omega\mathcal S$. More precisely, 
\begin{equation}
\label{eq:Gzeta}
\frac{1}{L(u^\eps_z)}\ms J_\eps(u^\eps_z) = \partial_{z_1} (E_{k+\eps K}(u^\eps_z))\ms e_1+
\big(\displaystyle\fintm |\omega|^2\ms dx\big)^{-1}\ms \partial_{z_2} (E_{k+\eps K}(u^\eps_z))\ms \omega\ms;
\end{equation}
\item[$iv)$] As $\eps\to 0$, we have
\begin{equation}
\label{eq:zero_order}
E_{k+\eps K}(u^\eps_z)- E_{k+\eps K}(\omega_z)=o(\eps)
\end{equation}
uniformly on $\Omega$, together with the derivatives with respect to the variable $z$.
\end{itemize}
\end{Lemma}

\proof
In order to shorten formulae, for $r>0$, $m\in\{0,2\}$ and $\delta>0$
we write
\begin{gather*}
\Omega_r=\{z\in\R^2~|~\text{dist}(z,\Omega)<r\}\ms,\\
C^m=C^m(\S^1,\R^2)\ms,\quad
\mathcal U_\delta:=\{\eta\in C^2~|~|\eta(x)|<\delta~~\text{for any $x\in\S^1$}~\}~\!.
\end{gather*}

Take $r, \delta>0$  small enough, so that $\overline\Omega_{2r}\subset\H^2$ and 
$\omega_z+\eta\in C^2(\S^1,\H^2)\setminus\H^2$
for any $z\in\overline\Omega_{2r}$, $\eta\in \mathcal U_\delta$.
Consider the differentiable function 
$$
\mathcal F:(\R\times\Omega_{2r})\times\mathcal U_\delta\!\times\!(\R\!\times\!\R^2)\to C^0\!\times\! (\R\!\times\!\R^2)~,\quad
\mathcal F=\big(\mathcal F_1,\mathcal F_2),
$$
whose coordinates
$$
\mathcal F_1:(\R\times\Omega_{2r})\times\mathcal U_\delta\!\times\!(\R\!\times\!\R^2)\to C^0\ms,\quad
\mathcal F_2:(\R\times\Omega_{2r})\times\mathcal U_\delta\!\times\!(\R\!\times\!\R^2)\to\R\!\times\!\R^2
$$
are given by 
\begin{eqnarray*}
\mathcal F_1(\eps,z;\eta;t,\vartheta)&=& J_\eps(\omega_z+\eta)-t\omega'-\vartheta_1e_1-\vartheta_2\omega,\\
\mathcal F_2(\eps,z;\eta;t,\vartheta)&=&\Big(\fintm\eta\cdot\omega'~\!dx,\fintm\eta_1~\!dx,\fintm\eta\cdot\omega~\!dx\Big).
\end{eqnarray*}
Take $z\in\Omega_{2r}$ and notice that 
$\mathcal F(0,z;0;0,0)=0$ because $J_0(\omega_z)=0$. 
The next goal is to solve the equation $\mathcal F(\eps,z;\eta;t,\vartheta)=(0,0)$ 
in a neighborhood of $(\eps,z)=(0,z)$, $(\eta;t,\vartheta)=(0;0,0)$ via the implicit function theorem.
Let
$$\mathcal L=(\mathcal L_1,\mathcal L_2):C^2\times(\R\times\R^2)\to C^0\times(\R\times\R^2)$$
be the differential of $\mathcal F(0,z;\ms\cdot\ms;\ms\cdot\ms,\ms\cdot\ms)$ computed at 
$(\eta;t,\vartheta)=(0;0,0)\in C^2\times(\R\!\times\!\R^2)$. We need to prove that $\mathcal L$ is invertible.
Explicitly, we have 
$$
\begin{array}{ll}
\mathcal L_1:C^2\times(\R\times\R^2)\to C^0, &\mathcal L_1(\f;a,{p})= J_0'(\omega_z)\f-a\omega'-{p}_1e_1-{p}_2\omega\\
&\\
\mathcal L_2:C^2\times(\R\times\R^2)\to\R\times\R^2, &\mathcal L_2(\f;a,{p})=
\displaystyle{\Big(\fintm\f\cdot\omega'~\!dx,\fintm\f_1~\!dx,\fintm\f\cdot\omega~\!dx\Big)}.
\end{array}
$$ 
If
$\mathcal L_1(\f;a,{p})=0$ then $J_0'(\omega_z)\f\in T_\omega\mathcal S$, hence 
$\f\in T_\omega\mathcal S$ by $ii)$ in Lemma \ref{L:nondegen}. If $\mathcal L_2(\f;a,{p})=0$
then $\f\in T_\omega\mathcal S^\perp$. Therefore, the operator $\mathcal L$ is 
injective.

\medskip

To prove surjectivity take  $u\in C^0, (b,q)\in\R\times\R^2$. We have to find 
$\f\in C^2$, $(a,p)\in\R\times\R^2$ satisfying $\mathcal L_1(\f;a,{p})=u$ and $\mathcal L_2(\f;a,{p})=(b,q_1,q_2)$, that is,
\begin{gather}
\label{eq:su1}
 J_0'(\omega_z)\f=u+a\omega'+{p}_1e_1+{p}_2\omega\\
\label{eq:su2}
\fintm\f\cdot\omega'~\!dx=b~,\quad \fintm\f_1~\!dx={q}_1~,\quad \fintm\f\cdot\omega~\!dx={q}_2.
\end{gather}
By (\ref{eq:autoaggiunto}),  for any $\f\in C^2$, $\tau\in T_\omega \mathcal S=\langle\omega',e_1,\omega\rangle=\ker\!J_0'(\omega_z)$ we have
$$
\fintm J_0'(\omega_z)\f\cdot\tau~\!dx=\fintm J_0'(\omega_z)\tau\cdot\f~\!dx=0.
$$
Thus the unknowns $a\in\R$ and ${p}=({p}_1,{p}_2)\in\R^2$ are determined by the condition
\begin{equation}
\label{eq:u_perp}
\fintm u\cdot\tau\ms dx+a\fintm \omega'\cdot\tau\ms dx+{p}_1\fintm e_1\cdot\tau\ms dx+{p}_2\fintm\omega\cdot\tau\ms dx=0
\quad\text{for any $\tau\in T_\omega \mathcal S$.}
\end{equation}
Now we look for the $L^2$ projection of the unknown function $\f$ on $T_\omega\mathcal S$ and its $L^2$ projection
on $T_\omega\mathcal S^\perp$. The tangential component
$\f^\top\in T_\omega\mathcal S=\langle \omega', e_1,\omega\rangle$  is uniquely determined by (\ref{eq:su2}). 
Next, we 
notice that 
$u+ a\omega'+{p}_1 e_1+{p}_2\omega\in   T_\omega \mathcal S^\perp$ by (\ref{eq:u_perp});
then we use $iii)$ in Lemma \ref{L:nondegen} to find $\f^\perp \in  C^2\cap T_\omega\mathcal S^\perp$ 
such that 
$$
 J_0'(\omega_z)\f^\perp=u+a\omega'+{p}_1e_1+{p}_2\omega~\!.
$$
The function $\f=\f^\top+\f^\perp$ solves (\ref{eq:su1}) because 
$J_0'(\omega_z)\f=J_0'(\omega_z)\f^\perp$, and surjectivity is proved.

We can now apply the implicit function theorem
for any fixed $z\in \Omega_{2r}$. Actually, thanks
a compactness argument we have that there exist
$\eps'>0$ and (uniquely determined) $C^1$ functions
$$
\begin{array}{lll}
\eta:(-\eps',\eps')\times\Omega_r\to \mathcal U_\delta\subset C^2~~&t:(-\eps',\eps')\times\Omega_r\to \R~~
&\vartheta:(-\eps',\eps')\times\Omega_r\to \R^2\\
\eta:(\eps,z)\mapsto \eta^\eps(z)~~
&t:(\eps,z)\mapsto t^\eps(z),~
&\vartheta:(\eps,z)\mapsto \vartheta^\eps(z)
\end{array}
$$
such that 
$$\eta^0(z)=0~,\quad t^0(z)=0~,\quad\vartheta^0(z)=0~,\quad
\mathcal F(\eps,z;\eta^\eps(z);t^\eps(z),\vartheta^\eps(z))=0.
$$
We introduce the $C^1$ function
$$
(-\eps',\eps')\times\Omega_r\to C^2(\S^1,\H^2)\setminus\H^2~,\quad (\eps,z)\mapsto u^\eps_z:=\omega_z+\eta^\eps(z)~\!,$$ 
that clearly satisfies $u^0_z=\omega_z$. Since $\omega_z$ is embedded, then $u^\eps_z$ is embedded as well,
provided that $\eps'$ is small enough. Moreover we have
\begin{gather}
\label{eq:natural1}
 J_\eps(u^\eps_z)=t^\eps(z)\omega'+\vartheta_1^\eps(z)e_1+\vartheta_2^\eps(z)\omega\in T_\omega\mathcal S\\
\label{eq:natural2}
\fintm(u^\eps_z-\omega_z)\cdot\omega'~\!dx=\fintm(u^\eps_z-\omega_z)\cdot e_1~\!dx=\fintm(u^\eps_z-\omega_z)\cdot\omega~\!dx=0,
\end{gather}
and (\ref{eq:natural2}) shows that $ii)$ is fulfilled.

Since  integration by parts gives
$$
\fintm\omega_z\cdot\omega'\ms dx=0~,\quad \fintm\omega_z\cdot e_1\ms dx=z_1~,\quad 
\fintm\omega_z\cdot\omega\ms dx=z_2\fintm|\omega|^2\ms dx\ms,
$$
we can rewrite the orthogonality conditions  (\ref{eq:natural2})
in the following, equivalent way:
\begin{equation}
\label{eq:natural3}
\fintm u^\eps_z\cdot\omega'~\!dx=0~,\quad \fintm u^\eps_z\cdot e_1~\!dx= z_1~,\quad \fintm u^\eps_z\cdot\omega~\!dx=z_2\fintm|\omega|^2~\!dx.
\end{equation}

Our next aim is to show that $t^\eps(z)=0$ for any $z\in\overline\Omega$, provided that $\eps$ is small enough.
We have that $\|(u^\eps_z)'-\omega'_z\|_\infty=o(1)$ as $\eps\to 0$, uniformly for $z\in\overline\Omega$.
Thus
$$
\fintm (u^{\eps}_{z})'\cdot\omega'~\!dx=\fintm \omega'_{z}\cdot\omega'~\!dx+o(1)=z_2 \fintm |\omega'|^2~\!dx+o(1).
$$
In particular, there exists $\beps\in(0,\eps')$ such that
$\displaystyle{\int_{\S^1} (u^{\eps}_{z})'\cdot\omega'~\!dx}$ is bounded away from $0$ if $(\eps,z)\in[-\beps,\beps]\times
\overline\Omega$.
On the other hand, using (\ref{eq:eps_invariance}), (\ref{eq:natural1}), integration by parts and (\ref{eq:natural3}), we 
have
\begin{eqnarray*}
0&=& \fintm J_\eps(u^{\eps}_{z})\cdot (u^{\eps}_{z})'~\!dx\\
&=&
t^\eps(z)\fintm (u^{\eps}_{z})'\cdot\omega'~\!dx+\vartheta_1^\eps(z)\fintm (u^{\eps}_{z})'\cdot e_1~\!dx+
\vartheta_2^\eps(z)\fintm (u^{\eps}_{z})'\cdot\omega~\!dx\\
&=&
t^\eps(z)\fintm (u^{\eps}_{z})'\cdot\omega'~\!dx-
\vartheta_2^\eps(z)\fintm u^{\eps}_{z}\cdot\omega'~\!dx= t^\eps(z)\fintm (u^{\eps}_{z})'\cdot\omega'~\!dx.
\end{eqnarray*}
We  see that $t^\eps(z)=0$ for any $(\eps,z)\in[-\beps,\beps]\times\overline\Omega$, and therefore
\begin{equation}
\label{eq:natural4}
  J_\eps(u^\eps_z)=\vartheta_1^\eps(z)e_1+\vartheta_2^\eps(z)\omega\ms.
\end{equation}
Now we compute the derivatives of the function $z\mapsto E_{k+\eps K}(u^\eps_z)$
via (\ref{eq:dL_eps}) and (\ref{eq:natural4}). For $j=1,2$ we obtain
\begin{eqnarray*}
 L(u^\eps_z)\partial_{z_j} (E_{k+\eps K}(u^\eps_z))
 &=&L(u^\eps_z)E'_{k+\eps K}(u^\eps_z)\partial_{z_j}u^\eps_z=
\fintm J_\eps(u^\eps_z)\cdot \partial_{z_j}u^\eps_z\ms dx\\
&=&
\vartheta_1^\eps(z)\fintm \partial_{z_j}u^\eps_z\cdot e_1\ms dx+ \vartheta_2^\eps(z)\fintm \partial_{z_j}u^\eps_z\cdot\omega\ms dx\\
&=&\vartheta_1^\eps(z)\partial_{z_j} \Big(\fintm u^\eps_z\cdot e_1\ms dx\Big)+
 \vartheta_2^\eps(z)\partial_{z_j}\Big(\fintm u^\eps_z\cdot\omega\ms dx\Big)\ms.
\end{eqnarray*}
Then we use (\ref{eq:natural3}) to infer
$$L(u^\eps_z)\partial_{z_1} (E_{k+\eps K}(u^\eps_z))=\vartheta_1^\eps(z)\ms,
\quad
L(u^\eps_z)\partial_{z_2} (E_{k+\eps K}(u^\eps_z))=\vartheta_2^\eps(z)\big(\fintm|\omega|^2\ms dx\big)\ms,
$$
that compared with (\ref{eq:natural4}) give 
(\ref{eq:Gzeta}).

\medskip

It remains to prove $iv)$. Take $z\in\overline\Omega$ and consider the function
$$
f_z(\eps)=E_{k+\eps K}(u^\eps_z)= E_k(u^\eps_z)+\eps A_K(u^\eps_z)\ms,
\quad f_z\in C^1(-\beps,\beps).
$$
Clearly $f_z(0)=E_k(\omega_z)$. To compute $f'_z(0)$ notice that
$\partial_\eps u^\eps_z$
remains bounded in $C^2(\overline\Omega)$ as $\eps\to 0$, because
the function $(\eps,z) \mapsto u^\eps_z$ is of class $C^1$. Thus
$A'_K(u^\eps_z)(\partial_\eps u^\eps_z)$ remains bounded as well. Further,
$E'_k(u^\eps_z)\to E'_k(\omega_z)=0$ in the norm operator because $u^\eps_z\to\omega_z$ in $C^2$
and since $\omega_z$ is a $k$-loop. 
We infer that
$$
f'_z(0)=E'_k(\omega_z)(\partial_\eps u^\eps_z)+A_K(u^\eps_z)+o(1)= A_K(\omega_z)+o(1)
$$ 
uniformly on $\overline\Omega$. In fact we proved that
$$f_z(\eps)=E_{k+\eps K}(u^\eps_z)=E_k(\omega_z)+\eps A_K(u^\eps_z)+o(1)$$
uniformly on $\overline\Omega$ as $\eps\to0$. That is, (\ref{eq:zero_order}) holds true  ''at the zero order''. 

\medskip

To conclude the proof we have to handle
$\partial_{z_j}\big(E_{k+\eps K}(u^\eps_z)- E_{k+\eps K}(\omega_z)\big)$ for $j=1,2$.
Since $J_\eps(u)=J_0(u)+\eps L(u)u_2^{-2}K(u)~\!iu'$, we can rewrite (\ref{eq:Gzeta}) as follows, 
\begin{align}
\nonumber
\partial_{z_1} (E_{k+\eps K}(u^\eps_z))\ms e_1+
\big(\displaystyle\fintm &|\omega|^2\ms dx\big)^{-1}\ms \partial_{z_2} (E_{k+\eps K}(u^\eps_z))\ms \omega\\
\label{eq:multline}
=&\ms \frac{1}{L(u^\eps_z)}J_0(u^\eps_z)+\eps (u^\eps_z)_2^{-2}K(u^\eps_z)i(u^\eps_z)'.
\end{align}
Recall that $J_0(u^\eps_z)$ is orthogonal to $e_1$ in $L^2$, see the second identity in
(\ref{eq:inv0}). We test (\ref{eq:multline}) with $e_1$ to obtain
\begin{equation}
\label{eq:domenica}
\partial_{z_1}\big(E_{k+\eps K}(u^\eps_z)\big)=
\eps\fintm  (u^\eps_z)_2^{-2}K(u^\eps_z)e_1\cdot i(u^\eps_z)'\ms dx=
\eps A'_K(u^\eps_z)e_1
\end{equation}
by (\ref{eq:dAK}). Since $\partial_{z_1}\big(E_{k+\eps K}(\omega_z)\big)=\partial_{z_1}\big(E_k(\omega)+\eps A_K(\omega_z)\big)
=\eps A'_K(\omega_z)e_1
$, we get
$$
\partial_{z_1}\big(E_{k+\eps K}(u^\eps_z)- E_{k+\eps K}(\omega_z)\big)=
\eps\big( A'_K(u^\eps_z)e_1 -  A'_K(\omega_z)e_1\big)=o(\eps)
$$
because of the continuity of $A_K'(\cdot)$ and since $u^\eps_z\to \omega_z$.

To handle the derivative with respect to $z_2$ we test (\ref{eq:multline}) with $u^\eps_z$.
Since $J_0(u^\eps_z)$ is orthogonal to $u^\eps_z$  in $L^2$  by (\ref{eq:inv0}), using also
(\ref{eq:natural3}) we obtain
$$
z_1\partial_{z_1}\big(E_{k+\eps K}(u^\eps_z)\big)+z_2 \partial_{z_2}\big(E_{k+\eps K}(u^\eps_z)\big)
=
\eps\fintm  (u^\eps_z)_2^{-2}K(u^\eps_z)u^\eps_z\cdot i(u^\eps_z)'\ms dx=
\eps A'_K(u^\eps_z)u^\eps_z~\!,
$$
that compared with (\ref{eq:domenica}) gives
$$
z_2 \partial_{z_2}\big(E_{k+\eps K}(u^\eps_z)\big)= 
\eps A'_K(u^\eps_z)(u^\eps_z-z_1e_1)\ms.
$$
From
$z_2\partial_{z_2}\big(E_{k+\eps K}(\omega_z)\big)=z_2\partial_{z_2}\big(E_k(\omega)+\eps A_K(\omega_z)\big)
=z_2\eps A'_K(\omega_z)\omega=\eps A'_K(\omega_z)(\omega_z-z_1e_1),
$
we conclude that 
$$
z_2\partial_{z_2}\big(E_{k+\eps K}(u^\eps_z)- E_{k+\eps K}(\omega_z)\big)
=
\eps\big( A'_K(u^\eps_z)(u^\eps_z-z_1e_1) -  A'_K(\omega_z)(\omega_z-z_1e_1)\big)
=o(\eps)\ms.
$$
The lemma is completely proved.
\QED

\subsection{Existence results}

\noindent{\bf Proof of Theorem \ref{T:main}.}
We are assuming that there exists $r>0$ such that any function ${G}\in C^1(\overline A)$
satisfying 
$\|{G}+{F^{\!K}_k}\|_{C^1(\overline A)}<r$  has a critical point in $A$. We recall also formula
(\ref{eq:pioggia}), that in particular gives
\begin{equation}
\label{eq:Fkk}
E_{k+\eps K}(\omega_z)=E_{k}(\omega_z)+\eps A_{K}(\omega_z)
=E_{k}(\omega)-\frac{\eps}{2\pi}\ms F^{\!K}_k(z)\ms.
\end{equation}

Take an open set $\Omega\in \H^2$ such that $A\Subset \Omega\Subset\H^2$, and
let $(\eps,z)\mapsto u^\eps_z$, $(\eps,z)\in[-\beps,\beps]\times\overline\Omega$ 
be the function given by Lemma \ref{L:reduction}. 
For $\eps\neq 0$ consider the function
$$
G^\eps(z)= \frac{2\pi}{\eps}(E_{k+\eps K}(u^\eps_z)-E_{k}(\omega))
$$
and use (\ref{eq:Fkk}) together with $iv)$ in Lemma \ref{L:reduction} to get
$$
\|G^\eps+F^{\!K}_k\|_{C^1(\overline A)}=
\frac{2\pi}{|\eps|}\big\|E_{k+\eps K}(u^\eps_z)-E_{k+\eps K}(\omega_z)\big\|=o(1)
$$
as $\eps\to 0$. 
We see that for $\eps$ small enough the function $G^\eps$
has a critical point $z^\eps\in A$. Since the derivatives of the function $z\mapsto E_{k+\eps K}(u^\eps_z)$
vanish at $z=z^\eps$, then  $J_\eps(u^\eps_{z^\eps})=0$ by (\ref{eq:Gzeta}). That is,
$u^\eps_{z^\eps}$ is
and embedded $k+\eps K$ loop.

The last conclusion in Theorem \ref{T:main} follows via a simple compactness argument
and thanks to Theorem \ref{T:necessary}. 
\QED

In the next result we apply Theorem \ref{T:main} to obtain the
existence of $k+\eps K$-loops that shrink to a stable critical point for the curvature function $K$, as $k\to\infty$.

\begin{Theorem}
\label{T:main2}
Let $K\in C^1(\H^2)$. Assume  that $K$ has a stable 
critical point in an open set $A\Subset \H^2$. 
There exists $k_0>1$ such that for any  $k>k_0$ 
and for every $\eps$ close enough to $0$, {there exists an embedded  $(k+\eps K)$-loop.}

Moreover, {let $k_h\to\infty, \eps_h\to 0$ be given sequences. There exist
subsequences $k_{h_j}, \eps_{h_j}$, a point $z_\infty\in\overline A$
that is critical for $K$,  and an embedded  $(k_{h_j}+\eps_{h_j} K)$-loop 
$u^{j}$ 
such that $u^j$ converges in $C^2(\S^1,\H^2)$ to the constant curve $z_\infty$,
as $h\to \infty$.}
\end{Theorem}

\proof
Recall that $R_k=({k^2-1})^{-1/2}$.
In order to simplify notations we put
$$
z^k:=(z_1, {kR_k}z_2)=z+(kR_k-1)z_2 e_2\quad \text{for $z=(z_1,z_2)\in\H^2$.}
$$
Since $D^\H_{\rho_k}(z)= D_{R_k z_2}(z^k)$ we have 
\begin{equation}
\label{eq:pioggia2}
F^{\!K}_k(z)=\int\limits_{D_{R_kz_2}(z^k)} p_2^{-2}K({p})~\!dp=\int\limits_{D_{R_k}(0)} 
(q_2+kR_k)^{-2}K({z_2q+z^k})~\!dq~\!.
\end{equation}
We put
$\phi_K(q)=q_2^{-2}K(q)$
and rewrite (\ref{eq:pioggia2}) as follows:
$$
\frac{1}{\pi R_k^2 z_2^2}~\!F^{\!K}_k(z)=\fint\limits_{D_{R_k}(0)} \phi_K(z_2q+z^k)~\!dq\ms.
$$
Trivially $kR_k=k/\sqrt{k^2-1}\to 1$ and  $|z^k-z|=(kR_k-1)z_2\to 0$ uniformly on $\overline A$,  as $k\to \infty$.
Since $\phi_K\in C^1(\H^2)$, it is easy to show that
$$\frac{1}{\pi R_k^2 z_2^2}~\!F^{\!K}_k(z)\to \phi_K(z)=\frac{1}{z_2^2}K(z)$$ in $C^1(\overline A)$.
It follows that for $k$ large enough,
$F^{\!K}_k$ has stable critical point in $A\Subset \H^2$. 
Theorem \ref{T:main} applies and gives the conclusion of the proof.
\QED

\appendix

\section{Loops in the Euclidean plane}
    
The  argument we used to prove Theorem \ref{T:main} applies also
in the easier Euclidean case. 
It is well known that the only embedded loops in $\R^2$ having prescribed
constant curvature $k>0$ are circles of radius $1/k$. 
We take as a reference
circle the loop
$$
\omega(x)=\frac{1}{k} \ms x~,\qquad x\in\S^1\subset\R^2,
$$
that solves 
$$
u''=L(u) k\ms iu'~,\quad\text{where}\quad 
L(u):=\Big(\fintm |u'|^2\ms dx\Big)^\frac12
$$
(in fact, $L(\omega)k=1$ and $\omega''=-\omega=i\omega'$).

Let $K\in C^1(\R^2)$ be given. If a nonconstant function $u\in C^2(\S^1,\R^2)$ solves
\begin{equation}
\label{eq:Euclideo}
u''=L(u)(k+\eps K(u))\ms iu'\ms,
\end{equation}
then $|u'|=L(u)$ is constant, and $u$ parameterizes a loop in $\R^2$
having Euclidean curvature $k+\eps K$ at each point. Further, problem (\ref{eq:Euclideo}) admits
a variational structure, see \cite{BCG}, \cite{MuLoop}. More precisely, its nonconstant solutions are critical points
of the energy functional
$$
E_{k+\eps K}(u)=\Big(\fintm |u'|^2\ms dx\Big)^\frac12+\eps \fintm Q(u)\cdot iu'~,\quad u\in C^2(\S^1,\R^2)\setminus\R^2,
$$
where the vectorfield $Q\in C^1(\R^2,\R^2)$ satisfies $\div Q=K$. 

Arguing as for Theorem \ref{T:necessary} one can prove a necessary conditions for the
existence of solutions to (\ref{eq:Euclideo}) for $\eps=\eps_h\to 0$. 

\begin{Theorem}
\label{T:necessaryE}
Let $u_h$ be a $(k+\eps_hK)$-loop solving 
(\ref{eq:Euclideo}) for $\eps=\eps_h$, and assume that

\medskip

\centerline{$ L(u_h)\to L_\infty>0$,
\qquad $u_h\to U$ uniformly, for some $U\in C^0(\S^1,\R^2)$.}

\medskip

\noindent
Then $U(x)=\omega\big(\xi x^\mu)+z$ for some $\mu\in\mathbb N$, $\xi\in\S^1$ and $z\in\R^2$, 
that is a critical point for the Melnikov function
$$
F^{\!K}_k(z)=\int\limits_{D_{\frac1k}(z)}K(q)\ms dq\ms,\quad F^{\!K}_k:\R^2\to\R\ms.
$$
\end{Theorem}

In the Euclidean case we have the following existence result.

\begin{Theorem}
\label{T:mainE}
Let $k>0$ and $K\in C^1(\R^2)$ be given. Assume  that $F^{\!K}_k$ has a stable 
critical point in an open set $A\Subset \R^2$. 
Then for every $\eps\in\R$ close enough to $0$, there exists an embedded $(k+\eps K)$-loop $u^\eps:\S^1\to\R^2$.

Moreover, any sequence $\eps_h\to 0$ has a subsequence $\eps_{h_j}$ such that 
$u^{\eps_{h_j}}\to \omega_{z_0}$ in $C^2(\S^1,\R^2)$ as $j\to\infty$, where
$z_0\in A$ is a critical point for $F^{\!K}_k$.
\end{Theorem} 

\noindent
{\bf Sketch of the proof.}
We introduce the  $3$-dimensional space of embedded solutions to the
unperturbed problem, namely
$$
\mathcal S=\big\{ \omega\circ\xi+z~|~\xi\in\S^1~,~~z\in\R^2~\big\},
$$
and the functions $J_\eps: C^2(\R,\R^2)\setminus\R^2\to C^0(\R,\R^2)$, $\eps\in\R$,
given by
$$
J_\eps(u)=-u'' +L(u)(k+\eps K(u))~\! iu'= J_0(u)+L(u) K(u)~\! iu'.
$$
We have $\mathcal S\subset \{J_0=0\}$. Since 
$\displaystyle{J'_0(\omega+z)\f= -\f''+\ms i\f'-k^2\big(\fintm \f\cdot\omega\ms dx\big)\omega}$, 
it is quite easy to check that
$$
T_{\omega+z}\mathcal S=\langle \omega', e_1, e_2\rangle=\ker J'_0(\omega+z),
$$
and that $J'_0(\omega+z):T_{\omega+z}\mathcal S^\perp\to T_{\omega+z}\mathcal S^\perp$
is invertible. The remaining part of the proof runs with minor changes. 
\QED

Theorem \ref{T:main2} has its Euclidean correspondent as well. We omit the proof of the next result.

\begin{Theorem}
\label{T:main2E}
Let $K\in C^1(\R^2)$. Assume  that $K$ has a stable 
critical point in an open set $A\Subset \R^2$. 
Then there exists $k_0>1$ such that for any fixed $k>k_0$, 
and for every $\eps$ close enough to $0$, there exists an embedded $(k+\eps K)$-loop $u^{k,\eps}:\S^1\to\R^2$.

Moreover, there exist sequences $k_h\to\infty$, $\eps_h\to 0$ such that
$u^{k_h,\eps_{h_j}}\to \omega_{z_0}$ in $C^2(\S^1,\R^2)$ as $j\to\infty$, where
$z_0\in A$ is a critical point for $K$.
\end{Theorem}

\footnotesize
\label{References}

\end{document}